\documentclass[titlepage,11pt]{article}
\usepackage{latexsym, tikz}
\usepackage{amsfonts}
\usepackage{amsmath}
\usepackage{textcomp}
\usetikzlibrary{graphs}
\usepackage{float}

\newcommand\blackslug{\hbox{\hskip 1pt \vrule width 4pt height 8pt depth 1.5pt
        \hskip 1pt}}
\newcommand\bbox{\hfill \quad \blackslug \bigbreak}
\def\dd{\hbox{-}}
\def\cc{\hbox{-}\cdots\hbox{-}}
\def\ll{,\ldots,}

\title{Finding a shortest odd hole}
\author{Maria Chudnovsky\thanks{This material is based upon work supported in part by the U. S. Army
Research Office under grant   number W911NF-16-1-0404, and supported by  NSF grant DMS 1763817.}\\
Princeton University, Princeton, NJ 08544
\\
\\
Alex Scott\\
Mathematical Institute, University of Oxford, Oxford OX2 6GG, UK
\\
\\
Paul Seymour\thanks{Supported by AFOSR grant
A9550-19-1-0187, and NSF grant DMS-1800053.}\\
Princeton University, Princeton, NJ 08544}

\date{December 22, 2019; revised \today}

\newtheorem{thm}{}[section]

\newcommand{\Proof}{\noindent{\bf Proof.}\ \ }

\begin{document}
\maketitle
\begin{abstract}
An odd hole in a graph is a induced cycle with odd length greater than 3.
In an earlier paper (with Sophie Spirkl), solving a longstanding open problem,
we gave a polynomial-time algorithm to test if a graph has an odd hole.   We subsequently showed that, for every $t$, there is a polynomial time algorithm to test whether a graph contains an odd hole of length at least $t$.  In this paper, we give an algorithm that finds a shortest odd hole, if one exists.
\end{abstract}

\section{Introduction}
All graphs in this paper are finite and have no loops or parallel edges. A {\em hole} of $G$ is an induced subgraph of $G$
that is a cycle
of length at least four, and an {\em odd hole} is a hole of odd length.  An {\em antihole} of $G$ is an induced subgraph
whose complement graph is a hole.

The class of graphs that have no odd holes and no odd antiholes hase been heavily studied since the 1960s.  
Indeed, the ``strong perfect graph conjecture'' of Claude Berge \cite{berge} stated that if a graph and its complement both have no odd holes, then its chromatic
number equals its clique number.  Berge's conjecture was open for many years, until it was proved by two of us, with Robertson and Thomas \cite{perfect},
in the early 2000s.  The corresponding algorithmic question, of finding a polynomial-time algorithm to test if a graph is perfect, was settled around the same time:
two of us, with Cornu\'ejols, Liu and Vu\v skovi\v c \cite{bergealg}, gave a polynomial-time algorithm to test if a graph
has an odd hole or odd antihole, and so test for perfection. 

Excluding both odd holes and odd antiholes has strong structural consequences.  However, if we just exclude odd holes, then the resulting class of graphs appears to be  (in some sense) much less well-structured.   It was only recently
that two of us \cite{oddchi} proved that if a graph has no odd holes then its chromatic number is bounded
by a function of its clique number, resolving an old conjecture of Gy\'arf\'as \cite{gyarfas}.  The complexity of recognizing graphs with no odd holes was also open for some time.  
While the algorithm of \cite{bergealg} could test for the presence of an odd
hole or antihole in polynomial time, we were unable to separate the test for odd holes from the test
for odd antiholes, and the complexity of testing for an odd hole remained open.   Indeed, there was reason to suspect that a polynomial time algorithm might not exist, as Bienstock \cite{b1,b2}  showed
that testing if a graph has an odd hole containing a given vertex is NP-complete.  Surprisingly, the problem was recently resolved in the positive: with  Sophie Spirkl \cite{oddholetest}, we gave a polynomial time algorithm to test if a graph has an odd hole.  

For graphs that do contain odd holes, it is natural to ask what we can determine about their lengths.   For example, what can be said about the shortest and longest odd holes in a graph?  

It is easy to see that  finding a {\em longest} odd hole is NP-hard, by reduction from Hamiltonian Path with a specified start and end vertex (the idea is to subdivide every edge once and add an odd path between the chosen start and end vertices $x$ and $y$: a long $xy$ path in the original graph corresponds to a long odd hole in the new graph).  On the positive side, in an earlier paper \cite{longoddholetest} extending the methods of \cite{oddholetest}, we gave for every constant $t$ a polynomial time algorithm to test whether a graph $G$ contains an odd hole of length at least $t$ (although  the running time is a polynomial in $|G|$ with degree $\Theta(t)$).

In this paper, we consider the problem of finding a {\em shortest} odd hole.   Building on our earlier work with Sophie Spirkl~\cite{oddholetest}, we give an algorithm to find a shortest odd hole, if there is one.
Thus, the main result of the paper is:
\begin{thm}\label{mainthm}
There is an algorithm with the following specifications:
\begin{description}
\item [Input:] A graph $G$.
\item [Output:] Determines whether $G$ has an odd hole, and if so finds the minimum length of an odd hole.
\item [Running time:] $O(|G|^{14})$.
\end{description}
\end{thm}
We have not tried hard to reduce the exponent 14 to something smaller, and this might be possible, with extra complications, 
as was done in the final section of~\cite{oddholetest}; but the current algorithm
is complicated already, and our first priority is keeping it as simple as we can.

The algorithm of~\cite{oddholetest} allows us to determine whether a graph has an odd hole.
Asking for a {\em shortest} odd hole adds significant additional difficulty.
The algorithm of~\cite{oddholetest} came in three parts:
\begin{itemize}
\item First we test whether $G$ contains a ``jewel'' or ``pyramid''; these are two kinds of induced subgraph that can easily 
be detected and if one is present, $G$ has an odd hole. Henceforth we can assume that $G$ contains no jewel or pyramid.
\item Now we generate a ``cleaning list'', a list of polynomially-many subsets of $V(G)$, such that if $G$ has an odd hole,
then for some shortest odd hole $C$, one of the sets ($X$ say) is disjoint from $V(C)$ and contains all ``$C$-major'' vertices.
These are the vertices not in $V(C)$ but with several neighbours in $C$. This works in graphs that have no pyramid or jewel.
\item Third, for each $X$ in the cleaning list, we test whether $G\setminus X$ has a shortest odd hole $C$ 
without $C$-major vertices. There is an easy algorithm for this, that works in graphs that have no pyramid or jewel.
\end{itemize}
How can this be modified to output the minimum length of an odd hole?
The test for pyramids and jewels
used in the first step is the main problem: it will detect a pyramid or jewel if there is one, and thereby find some odd hole,
but not necessarily the shortest. We have to replace this with something else, and then adjust the second and third steps accordingly.

\section{Pyramids and jewels}

Let us give some definitions before we go on.
Let $v_0\ll v_3\in V(G)$ be distinct, and for $i = 1,2,3$ let $P_i$ be an induced path of $G$ between $v_0$ and $v_i$, such that

\begin{itemize}
\item $P_1,P_2,P_3$ are pairwise vertex-disjoint except for $v_0$;
\item at least two of $P_1,P_2,P_3$ have length at least two;
\item $v_1,v_2,v_3$ are pairwise adjacent; and
\item for $1\le i<j\le 3$, the only edge between $V(P_i)\setminus \{v_0\}$ and $V(P_j)\setminus \{v_0\}$ is the edge $v_iv_j$.
\end{itemize}
We call $P_1, P_2, P_3$ the {\em constituent paths} of the pyramid.
\begin{figure}[H]
\centering

\begin{tikzpicture}[scale=0.8,auto=left]
\tikzstyle{every node}=[inner sep=1.5pt, fill=black,circle,draw]

\node (v0) at (0,0) {};
\node (v1) at (-1,-3.5) {};
\node (v2) at (0,-3) {};
\node (v3) at (1,-3.5) {};
\tikzstyle{every node}=[]
\draw (v1) node [below]           {\footnotesize$v_1$};
\draw (v2) node [below]           {\footnotesize$v_2$};
\draw (v3) node [below]           {\footnotesize$v_3$};
\draw (v0) node [above]           {\footnotesize$v_0$};

\draw (v1) -- (v2);
\draw (v1) -- (v3);
\draw (v2) -- (v3);
\draw[dashed] (v0) -- (v1);
\draw[dashed] (v0) -- (v2);
\draw[dashed] (v0) -- (v3);

\end{tikzpicture}

\caption{A pyramid. Throughout, dashed lines represent paths, of indeterminate length.} \label{fig:pyramid}
\end{figure}
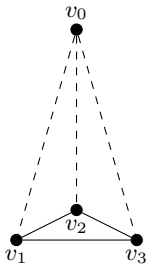

We call the subgraph induced on $V(P_1\cup P_2\cup P_3)$ a {\em pyramid}, with {\em apex} $v_0$ and {\em base} $\{v_1,v_2,v_3\}$.
If $G$ has a pyramid then $G$ has an odd hole (because two of the paths $P_1,P_2,P_3$ have the same length modulo two, and they induce
an odd hole).

If $X\subseteq V(G)$, we denote the subgraph of $G$ induced on $X$ by $G[X]$. If $X$ is a vertex or edge of $G$, or a set of vertices 
or a set of edges of $G$, we denote by $G\setminus X$ the graph obtained from $G$ by deleting $X$.
Thus, for instance, if $b_1b_2$ is an edge of a hole $C$, then $C\setminus \{b_1,b_2\}$ and $C\setminus b_1b_2$ are both paths,
but one contains $b_1,b_2$ and the other does not. If $P$ is a path, 
the {\em interior} of $P$ is the set of vertices of the path $P$ that are not ends of $P$.

We say that $G[V(P)\cup \{v_1,\ldots, v_5\}]$ is a {\em jewel} in $G$ 
if $v_1,\ldots, v_5$ are distinct vertices,
 $v_1v_2, v_2v_3,\allowbreak v_3v_4,\allowbreak v_4v_5, v_5v_1$ are edges, $v_1v_3,v_2v_4,v_1v_4$
are nonedges, and $P$ is a path of $G$ between $v_1,v_4$ such that $v_2,v_3,v_5$ have no neighbours
in the interior of $P$. 
(We do not specify whether $v_5$ is adjacent to $v_2,v_3$, but if it is adjacent to one and not the other, then $G$ also contains
a pyramid.)
\begin{figure}[H]
\centering

\begin{tikzpicture}[scale=0.8,auto=left]
\tikzstyle{every node}=[inner sep=1.5pt, fill=black,circle,draw]

\node (v1) at (-3,0) {};
\node (v2) at (-1,0) {};
\node (v3) at (1,0) {};
\node (v4) at (3,0) {};
\node (v5) at (0,2) {};
\tikzstyle{every node}=[]
\draw (v1) node [below left]           {\footnotesize$v_1$};
\draw (v2) node [below]           {\footnotesize$v_2$};
\draw (v3) node [below]           {\footnotesize$v_3$};
\draw (v4) node [below right]           {\footnotesize$v_4$};
\draw (v5) node [above]           {\footnotesize$v_5$};
\node  (P) at (0,-1.3)            {\footnotesize$P$};

\draw (v1) -- (v2) -- (v3)--(v4)--(v5)--(v1);
\draw[dotted] (v5) -- (v2);
\draw[dotted] (v5) -- (v3);
\draw[dashed] (v1) to [bend right = 60] (v4);

\end{tikzpicture}

\caption{A jewel. Throughout, dotted lines represent possible edges.} \label{fig:jewel}
\end{figure}
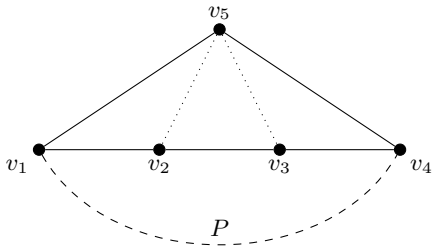

Every graph containing a pyramid or jewel has an odd hole,
and it was shown in~\cite{bergealg} that there is a polynomial-time algorithm to test if a graph contains a pyramid or jewel.
This was central to the algorithm of~\cite{oddholetest}, but it is no longer useful for us, as it stands.

But the test for jewels is easy to repair. Let us say an odd hole $C$ of $G$ is {\em jewelled} if either
\begin{itemize}
\item there is a four-vertex path of $C$ with vertices $c_1\dd c_2\dd c_4\dd c_5$ in order, and a vertex $c_3\in V(G)$ adjacent to 
$c_1$ and to $c_5$; or
\item there is a three-vertex path of $C$ with vertices $c_1\dd c_3\dd c_5$ in order, and two more vertices 
$c_2,c_4\in V(G)\setminus V(C)$, such that $c_1\dd c_2\dd c_4\dd c_3$ is an induced path.
\end{itemize}
There is a jewel in $G$ if and only if there is a jewelled odd hole in $G$.

\begin{thm}\label{testjewel}
There is an algorithm with the following specifications:
\begin{description}
\item [Input:] A graph $G$.
\item [Output:] Decides if there is a jewelled odd hole in $G$, and if so, finds a shortest one.
\item [Running time:] $O(|G|^7)$.
\end{description}
\end{thm}
\Proof We enumerate all five-tuples $(c_1\ll c_5)$ of vertices such that $c_1\dd c_2\dd c_4\dd c_5$ is an induced path and $c_3$
is adjacent to $c_1, c_5$. For each such choice we find a path $P$ of minimum length joining $c_1,c_5$ whose interior contains
no neighbours of $c_2,c_3$ or $c_4$, if there is one. If $P$ has odd length, we record the jewelled odd hole
$c_1\dd c_2\dd c_4\dd c_5\dd P\dd c_1$; and if $P$ is even, we record the jewelled odd hole $c_1\dd c_3\cc c_5\dd P\dd c_1$.
We output the shortest of all the recorded holes, or if there are none, report that no odd hole is jewelled. This proves \ref{testjewel}.~\bbox

(We could make this faster by adding complications to the algorithm, but there would be no gain in the 
overall running time of the main algorithm.) So for our shortest odd hole problem,
if some shortest odd hole happens to be jewelled, then the length of the shortest odd hole is the output of \ref{testjewel}.
This turns out to be good enough to replace the old test for jewels. Doing something similar for pyramids is a much greater challenge,
and is the main part of the paper.

\section{Handling pyramids}

If $G$ contains a pyramid, then some cycle of the pyramid is an odd hole of $G$. It turns out that we do not really need to
know that $G$ contains no pyramid; it is enough that there is no pyramid which includes a {\em shortest} odd hole of $G$. If we
could arrange that, then the remainder of the old algorithm could be used verbatim. Unfortunately we were unsuccessful.

But we can do something like it, which we will begin to explain in this section. If $P$ is a path, its {\em interior}
is the set of vertices of $P$ which have degree two in $P$, and is denoted by $P^*$.
Let $C$ be a shortest odd hole of $G$.
A vertex $v\in V(G)$ is {\em $C$-major} if there is no three-vertex path of $C$ containing all the neighbours of $v$ in $V(C)$
(and consequently $v\notin V(C)$); and $C$ is {\em clean} (in $G$) if no vertices of $G$ are $C$-major.
A $C$-major vertex is {\em big} if it has at least four neighbours in $V(C)$. It is easy to check that
\begin{thm}\label{notbig}
Let $C$ be a shortest odd hole of $G$. Let $v\in V(G)$ be $C$-major.
\begin{itemize}
\item If $v$ has at most three neighbours in $V(C)$, 
then $v$ has exactly three 
and exactly one pair of them are adjacent. 
\item If $v$ has exactly four neighbours in $V(C)$, then either exactly one pair of them are adjacent, or
$C$ is jewelled.
\end{itemize}
\end{thm}

Again, let $C$ be a shortest odd hole in $G$,
let $u,v \in V(C)$ be distinct and nonadjacent, and let $L_1,L_2$ be the two subpaths of $C$ joining $u,v$.
Suppose that there is a path $P$ of $G$, with ends $u,v$, such that 
\begin{itemize}
\item $|E(P)| <\min(|E(L_1)|,|E(L_2)|)$; and 
\item no big $C$-major vertex belongs to $V(P)$.
\end{itemize}
We call $P$ a {\em shortcut} for $C$. If $u,v\in V(G)$, $d_G(u,v)$ denotes the length of the shortest path of $G$ between $u,v$
(and $d_G(u,v)=\infty$ if there is no such path).

Let $H$ be a pyramid in $G$, with apex $a$ and base $\{b_1,b_2,b_3\}$, and constituent paths $P_1, P_2, P_3$, where $P_i$
is between $a, b_i$ for $i = 1,2,3$. Suppose that:
\begin{itemize}
\item the hole $a\dd P_1\dd b_1\dd b_2\dd P_2\dd a$ is a shortest odd hole $C$ of $G$; and
\item the length of $P_3$ is strictly less than the length of $P_i$ for $i = 1,2$.
\end{itemize}
In this case we call $H$ a {\em great pyramid} in $G$, and we call $V(P_3)\setminus \{a\}$ its {\em heart}.

The algorithm breaks into two parts: one will find a shortest odd hole if $G$ contains a great pyramid, and the other will find a shortest odd hole if $G$ does not contain a great pyramid.
More exactly, we will present two algorithms, as follows (a {\em 5-hole} means a hole of length five): the first is proved in sections
7--9, and the second in section 6.

\begin{thm}\label{findgreatpyr}
There is an algorithm with the following specifications:
\begin{description}
\item [Input:] A graph $G$.
\item [Output:] Outputs either an odd hole of $G$, or a statement of failure. If $G$ contains no 5-hole, and no jewelled shortest odd 
hole, and $G$ contains a great pyramid, 
the output will
be a shortest odd hole of $G$.
\item [Running time:] $O(|G|^{14})$.
\end{description}
\end{thm}

\begin{thm}\label{nogreatpyr}
There is an algorithm with the following specifications:
\begin{description}
\item [Input:] A graph $G$.
\item [Output:] Outputs either an odd hole of $G$, or a statement of failure. If $G$ contains no 5-hole, no jewelled shortest odd hole,
and no great pyramid, and $G$ contains an odd hole, the output will
be a shortest odd hole of $G$.
\item [Running time:] $O(|G|^9)$.
\end{description}
\end{thm}

The algorithm of \ref{findgreatpyr} was derived from
the algorithm in~\cite{bergealg} that tests if $G$ contains a pyramid; and that of \ref{nogreatpyr} 
is very similar to the main algorithm of~\cite{oddholetest}. 
We also need one more step:
\begin{thm}\label{find5hole}
There is an algorithm with the following specifications:
\begin{description}
\item [Input:] A graph $G$.
\item [Output:] Outputs either a 5-hole of $G$, or a statement of failure. If $G$ contains a 5-hole, 
the output will
be a shortest odd hole of $G$.
\item [Running time:] $O(|G|^5)$.
\end{description}
\end{thm}
\Proof We test all five-tuples of vertices of $G$.~\bbox

Let us derive the main result from these:

\bigskip

\noindent{\bf Proof of \ref{mainthm}.\ \ }
We have input a graph $G$. We run the algorithms of \ref{find5hole}, \ref{testjewel}, \ref{findgreatpyr} and \ref{nogreatpyr} on $G$, and for each
record the hole it outputs if there is one. If no hole is recorded, we report that $G$ has no odd hole; otherwise we output the shortest recorded hole.

To see correctness, there are five possibilities:
\begin{itemize}
\item $G$ has no odd hole: then each algorithm will report failure, and the output will be correct. 
\item $G$ has a 5-hole; then the output of \ref{find5hole} is a shortest odd hole.
\item $G$ has no 5-hole, and some shortest odd hole is jewelled: then \ref{testjewel} outputs a shortest odd hole.
\item $G$ has no 5-hole, and no jewelled shortest odd hole, and $G$ contains a great pyramid: then \ref{findgreatpyr} outputs a shortest odd hole.
\item $G$ has an odd hole, and has no 5-hole, no jewelled shortest odd hole, and no great pyramid: then \ref{nogreatpyr} outputs a shortest
odd hole.
\end{itemize}
In each of the last four cases, the recorded hole of smallest length is a shortest odd hole. This proves correctness of the algorithm.~\bbox

\section{Great pyramids}

In this section we prove three results very closely related to theorem 4.1 of~\cite{bergealg}:

\begin{thm}\label{shortcut}
Suppose that no shortest odd hole of $G$ is jewelled.
Let $P$ be a path with minimal interior such that for some shortest odd hole $C$ of $G$, $P$ is a shortcut for $C$.
Then the subgraph induced on $V(P\cup C)$ is a great pyramid in $G$ with heart $P^*$.
\end{thm}
\Proof Let $P$ have vertices $u\dd p_1\cc p_k\dd v$ in order, where $u,v\in C$, and $L_1,L_2$ are as in the definition of shortcut.
Since $|E(P)|<|E(L_1)|$, it follows that $k\ge 1$; and so
$\min(|E(L_1)|,|E(L_2)|\ge 3$. If $P_1$ is adjacent to $v$, then $p_1$ is $C$-major, but not big by hypothesis; and $k=1$ from the 
minimality of $P^*$; and the subgraph induced on $V(P\cup C)$ is a great pyramid, so the theorem holds.
Thus we may assume that $p_1$ is not adjacent to $v$, and hence $k\ge 2$. By the same argument, we may assume that no vertex of
$P^*$ is $C$-major.
Now $p_1$ may have more than one neighbour in $V(C)$, and the same for $p_k$, so let us choose $u,v$ to maximize 
$d_C(u,v)$. 

Assign
$C$ an orientation, clockwise say, and for any two distinct vertices $x,y$ in $C$, let $C(x,y)$ be the clockwise
path in $C$ from $x$ to $y$.
We may assume that $L_1 = C(u,v)$.
Let $C$ have vertices $c_1\cc c_{2n+1}$
in clockwise order, where $c_1 = u$ and $c_m = v$.
\\
\\
(1) {\em $P$ is an induced path of $G$ between $u,v$;  and no vertex of $P^*$ belongs to $V(C)$.}
\\
\\
The first claim is immediate from the minimality of $P^*$. For the second, suppose that $p_i\in V(C)$ say, where $1\le i\le n$.
From the symmetry we may assume that $p_i\in L_1^*$. From the minimality of $P^*$, the path $u\dd p_1\cc p_i$ is not a shortcut 
for $C$, and so its length is at least that one of $C(u,p_i), C(p_i,u)$. But $C(p_i,u)$ includes $L_2$ and so is longer than
$u\dd p_1\cc p_i$; and hence $C(u,p_i)$ has length at most $i$. Simularly $C(p_i,v)$ has length at most $k-i+1$; but then
$C(u,v)$ has length at most $k+1$, which is the length of $P$, a contradiction. This proves (1).

\bigskip

Since $C$ is odd, we may assume from the symmetry that $|E(L_1)|<|E(L_2)|$,
and therefore $m \le n+1$. 
From the hypothesis,
\[k+1 = d_G(u,v) < d_C(u,v) = m-1 \le n.\]
Since $k\ge 2$, it follows 
that $m\ge 5$ and so $n \ge 4$.
\\
\\
(2) {\em There are no edges between
$\{p_2\ll p_{k-1}\}$ and $C(v,u)^*$.}
\\
\\
For suppose not. Then for some $j$ with $m+1 \le j \le 2n+1$, there exist
paths $P_1,P_2$ from $c_j$ to $u,v$ respectively, with interior in $P^*$, both
strictly shorter than $P$. Suppose first that $j= 2n+1$.  Then
\[d_C(c_{2n+1},v) = \min (m, 2n+1-m) \ge m-1 >  |E(P)| > |E(P_2)|,\]
contrary to the minimality of $P^*$.
Thus $j\le 2n$ and similarly $j\ge m+2$. In particular, $P_1,P_2$ both have length at least two.
Now
\[|E(P_1)| + |E(P_2)| \le k+3 \le m \le 2n+2-m \le  (2n+2-j) + (j -m).\]
But $|E(P_1)|\ge 2n+2-j$ from the minimality of $P^*$, and similarly $|E(P_2)|\ge j-m$; so equality holds throughout.
In particular $k+3=m=n+1$, and $c_j$ is adjacent to $p_{2n+1-j}$ and to no other vertex in $P^*$. 
The lengths of $P,L_1$ differ by exactly one, and since $P\cup L_1$ is not an odd hole
(because it is shorter than $C$) it follows that some vertex $p_i$ of $P^*$ has a neighbour in $L_1^*$. 
Let $p_i$ be adjacent to $c_h$ where $2\le h\le m-1$. 
If $i>1$ then the path $c_h\dd p_i\cc p_k\dd c_m$ is not a shortcut over $C$, from the minimality of $P^*$, and
therefore $k-i+2\ge m-h =k+3-h$, and so $h-i\ge 1$. Similarly
(by exchanging $u,v$) it follows that if $i<k$ then $k-i\le m-h-1= k+2-h$, and so $h\le i+2$. 
The path $Q$ with vertices $c_h\dd p_i\cc p_{2n+1-j}\dd c_j$ has length $|2n+1-i-j|+2$.
We claim that $|E(Q)|<j-h$. 
If $i\le 2n+1-j$ then $i< k$ (since otherwise $p_k$ is adjacent to both $c_h, c_j$ and is therefore $C$-major, a contradiction),
and so 
$$|E(Q)|=|2n+1-i-j|+2=2n+3-i-j\le 2n+3-(h-2)-j=2n+5-h-j< j-h$$
since $j\ge m+2=n+3$. If $i>2n+1-j$, then 
$$|E(Q)|= |2n+1-i-j|+2= i+j+1-2n < j-h$$
because $h\le m-1=n$ and $i\le k=n-2$. Thus in either case $|E(Q)|<j-h$.
Similarly (by exchanging $u,v$)
it follows that $|E(Q)|<2n+1-h-j$, and therefore $Q$ is a shortcut for $C$,
contrary to the minimality of $P^*$. This proves (2).
\\
\\
(3) {\em Either $c_1$ is the only neighbour of $p_1$ in $C$, or $c_1,c_2$ are the only neighbours of
$p_1$ in $C$, or $m = n+1$ and $c_1,c_{2n+1}$ are the only neighbours of $p_1$ in $C$.
The analogous statement holds for $p_k$.}
\\
\\
For suppose first that $p_1$ has two nonadjacent neighbours $x,z \in V(C)$. Since $p_1$ is not $C$-major, we
may assume that $C(x,z)$ has length $2$ and contains all neighbours of $p_1$ in $C$.
Let $y$ be the middle vertex of $C(x,z)$; then $u\in \{x,y,z\}$, and since
$u,v$ are nonadjacent and $p_1,v$ are nonadjacent by (1), it follows that $v\ne x,y,z$. Now
$p_1\dd z \dd C(z,x) \dd x \dd p_1$ is a hole $C'$ of the same length as $C$, and hence is a shortest
odd hole. Suppose that $p_i$ is big $C'$-major for some $i\in \{2\ll k\}$. Since $p_i$ is not $C$-major, $p_i$
is adjacent to $p_1$, and so $i = 2$; and $p_i$ has three neighbours in $C$, since it has four in $C'$; and they are consecutive
since $p_2$ is not $C$-major, and so $C'$ is jewelled by \ref{notbig}, a contradiction.
This proves that $p_1$ does not have two nonadjacent neighbours in $V(C)$, and the same holds for $p_k$.

Since $p_1$ is adjacent to $c_1$, we may assume it is also
adjacent to $c_{2n+1}$, for otherwise the claim holds. Suppose that $m\le n$.
Then $d_C(c_{2n+1}, c_m)>d_C(c_1, c_m)$, contrary to the
choice of $u,v$ (maximizing $\min(|E(L_1)|,|E(L_2)|$).
Hence $m = n+1$. This proves (3).
\\
\\
(4) {\em We may assume that there are no edges between $P^*$ and $C(v,u)^*$.}
\\
\\
For suppose there are edges between $P^*$ and $C(v,u)^*$. From (2) and (3), we may assume that
$p_1$ is adjacent to $c_{2n+1}$ and $m = n+1$. Let $P'$ be the path $c_{2n+1}\dd p_1 \cc p_k\dd v$.
Then $P'$ is another shortcut for $C$, with the same interior as $P$.
Hence there is symmetry between $c_1,c_{2n+1}$, and from (2)
applied under this symmetry we deduce that 
there are no edges between $\{p_2\ll p_{k-1}\}$ and $C(u,v)^*$. Consequently there are no edges
between $P^*$ and $V(C)$ except for $p_1c_1,p_1c_{2n+1}$ and possibly edges incident with $p_k$.
By (3), and the symmetry between $c_1, c_{2n+1}$, we may assume that $p_k$ has no neighbours in $V(C)$
except $c_{n+1}$ and possibly $c_n$. Suppose that $p_k$ is adjacent to $c_n$. The holes
$$p_1\cc p_k\dd c_m\dd c_{m-1}\cc c_1\dd p_1$$
$$p_1\cc p_k\dd c_m\dd c_{m+1}\cc c_{2n+1}\dd p_1$$
are both shorter than $C$ and hence have even length; and so 
$(k+1)+(n-1)=k+n$ and $k+1+(2n+1-(n+1))= k+1+n$ are both even, which is impossible. 
Hence $c_{n+1}$ is the only neighbour of $p_k$ in $V(C)$; and so $G[V(C\cup P)]$ is a great pyramid.
This proves (4).
\\
\\
(5) {\em There are no edges between $\{p_2\ll p_{k-1}\}$ and $C(u,v)$.}
\\
\\
Suppose the claim is false;
then there exists $j$ with $2 \le j \le m-1$, and paths $P_1,P_2$ from $c_j$ to $u,v$
respectively, with interior in $P^*$, and both strictly shorter than $P$. Since
\[|E(P_1)| + |E(P_2)| \le k+3 \le m < j + (m-j+1)\]
it follows that either $P_1$ has length $<j$ or $P_2$ has length $< m-j+1$, and from the symmetry we may assume
the first. By the minimality of $P^*$, $P_1$ is not a shortcut, and so its length is exactly $j-1$. Thus
$c_j$ is adjacent to $p_{h}$, where $h=j-2$. By the same argument the only edge between $\{p_1\ll p_{h}\}$ and $\{c_j\ll c_m\}$
is the edge $p_{h}c_j$; and so by (4), the union of $P_1$ and the path $C(c_j,c_1)$ is a hole $C'$ say.
Thus $C'$ is a shortest odd hole. Suppose that $p_i$ is big $C'$-major, for some $i\in \{h+1\ll k\}$.
Then $p_i$ is adjacent to one of $p_1\ll p_h$, and therefore $i=h+1$; and $i$ has three neighbours in $V(C)$, consecutive;
and so it has exactly four in $V(C')$, and therefore $C'$ is jewelled, by \ref{notbig}, a contradiction.
From the minimality of $P^*$, it follows that $p_h\cc p_k\dd c_m$ is not a shortcut for $C'$, and so
$k-h+1\ge d_{C'}(p_h,c_m)=m-j+1$, that is, $k+2\ge m$ since $h=j-2$, a contradiction. This proves (5). 

\bigskip

Note that $c_{2n+1},c_1,c_2$ are all different from $c_{m-1},c_m,c_{m+1}$, since
$k\ge 2$.
From (2), (3) and (5) it follows that the only edges between $P^*$ and $V(C)$ are $p_1c_1$, $p_kc_m$,
possibly one edge from $p_1$ to one of $c_2,c_{2n+1}$, and possibly one edge from $p_k$
to one of $c_{m-1},c_{m+1}$.
If neither or both of the possible extra edges are present, there
is an odd hole shorter than $C$, a contradiction; so exactly one of the possible extra edges
is present. But then $G[V(C\cup P)]$ is a great pyramid. This proves
\ref{shortcut}.\bbox

\begin{thm}\label{reroute}
Let $G$ be a graph such that no shortest odd hole in $G$ is jewelled. 
Let $C$ be a shortest odd hole in $G$, and let $P$ be a path of $G$ with ends $u,v\in V(C)$, such that
$|E(P)|=d_C(u,v)$, and no vertex in $P^*$ is big $C$-major. Let $L_1, L_2$ be the paths of $C$ joining $u,v$, 
where $|E(L_1)|=|E(P)|$. Then either
\begin{itemize}
\item  $P\cup L_2$ is a shortest odd hole of $G$, or 
\item there is a great pyramid in $G$ with heart a proper subset of $P^*$, or
\item $G[V(P\cup C)]$ is a great pyramid with heart $P^*$.
\end{itemize}
\end{thm}
\Proof We proceed by induction on the length of $P$. The result is true if $P$ has length two, so we may assume that
$P$ has length at least three. We may assume that $G$ contains no great pyramid with heart a proper subset of $P^*$; so in particular,
for every shortest odd hole $C'$, every $C'$-major vertex in $P^*$ is big $C'$-major; and by \ref{shortcut} there is no shortcut 
for $C'$ with interior a proper subset of $P^*$.

In particular, $P$ is an induced path. Suppose
that $P\cup L_2$ is not a shortest odd hole. Then some vertex of $L_2^*$ is equal to or adjacent to some vertex of $P^*$.
Let $P$ have vertices $u\dd p_1\cc p_k\dd v$ in order, and let $C$ have vertices $c_1\dd c_2\cc c_{2n+1}\dd c_1$,
where $c_1=u$. We may assume that $c_{k+2}=v$ since $P,L_1$ have the same length; and so $k\le n-1$, since $L_1$ is shorter
than $L_2$. Choose $j$ with $k+3\le j\le 2n+1$ such that $c_j$ is equal or adjacent to some vertex $p_h$ in $P^*$. Thus there are
induced paths $P_1,P_2$ between $c_j$ and $u,v$ respectively, with interior in $P^*$, such that 
$P_1$ has length at most $h+1$ and $P_2$ has length at most $k+2-h$. Consequently
$|E(P_1)|+|E(P_2)|\le k+3$. 

Suppose that $P_1$ is a shortcut for $C$. Then $P_1^* = P^*$, and by \ref{shortcut}, it follows that 
$G[V(P_1\cup C)]$ is a great pyramid with heart $P^*$, and the theorem holds. So we may assume that $P_1$ is not a shortcut for $C$,
So the length of $P_1$ is at least $d_C(u,c_j)=\min(2n+2-j,j-1)$.
But the length of $P_1$ is at most that of $P$, and so at most $k+1$; and $j-1> k+1$ since $j>k+2$. Consequently
$P_1$ has length at least $2n+2-j$, and so $h+1\ge 2n+2-j$. Similarly, by exchanging $u,v$ we may assume that $P_2$ has 
length at least $j-k-2$, and so $k+2-h\ge j-k-2$.
Adding, we deduce that
$$k+3\ge 2n-k\ge 2(k+1)-k=k+2.$$
Hence, either $h+1=2n+2-j$, or $k+2-h=j-k-2$. From the symmetry between $u,v$ we may
assume the first holds with equality and the second holds with an error of at most $1$. So $h+1=2n+2-j$ and $k+2-h\le j-k-1$. 
Since $P_1$ has length at most $h+1$ and at least $2n+2-j$, it follows that $P_1$ has length exactly
$h+1=2n+2-j$; and so $P_1$ is an induced path with vertices $c_1\dd p_1\cc p_h\dd c_j$ in order. (In particular, $p_h\ne c_j$,
and $j\ne 2n+1$.)

Suppose that the length of $P_1$ is less than the length of $P$. No vertex of $P_1^*$ is $C$-major, so from the inductive
hypothesis, the union of $P_1$ and the path $c_1\dd c_2\cc c_j$ is a shortest odd hole $C'$ say. 
Suppose that one of $p_{h+1}\ll p_k$, say $p_i$, is big $C'$-major. Since $p_i$ is not big $C$-major, $p_i$ has a neighbour in 
$\{p_1\ll p_h\}$, and hence $i=h+1$. But $p_i$ has at least four neighbours in $V(C')$, and so it has at least three in 
$V(C\cap C')$. They are consecutive since $p_i$ is not $C$-major; but this contradicts \ref{notbig} applied to $C'$. This proves
that none of $p_{h+1}\ll p_k$ is big $C'$-major. But then the path $p_h\dd p_{h+1}\cc p_k\dd c_{k+2}$ is a shortcut for $C'$,
since $d_{C'}(p_h,c_{k+2})$ is the minimum of $j-k-1, h+k+1$, and both the latter are greater than $k-h+1$. This is
a contradiction, and so the length of $P_1$ equals that of $P$; that is, $h=k$. But $h+1=2n+2-j$ and so $j=2n+1-k$.
Since $p_k$ is adjacent to both $c_{k+2}$ and $c_j$, and $p_k$ is not $C$-major, it follows that $j\le k+4$, and so
$2n+1-k\le k+4$, that is, $n\le k+1$. But we already saw that $k\le n-1$, and so equality holds, and $j=2n+1-k=k+3$.

Thus we have proved so far that $k=n-1$, and $P^*$ is disjoint from $L_2^*$, and the only edges between $P^*$ and $L_2^*$
are either between $p_k,c_{k+3}$ or between $p_1,c_{2n+1}$, and we are assuming that at least one of these is present. 
If both are present then 
$p_1\cc p_k\dd c_{k+3}\dd c_{k+4}\cc c_{2n+1}\dd p_1$
is an odd hole shorter than $C$, a contradiction. So exactly one is present, say $p_1c_{2n+1}$. From the path 
$c_{2n+1}\dd p_1\cc p_k\dd c_{k+2}$ and the hole $C$, it follows that the edge $p_1c_1$ is the only edge between $P^*$ and
$\{c_1\ll c_{k+1}\}$; and so the subgraph induced on $V(C\cup P)$ is a great pyramid with heart $P^*$, and the theorem holds. 
This proves \ref{reroute}.~\bbox

\begin{thm}\label{majoreq}
Let $G$ be a graph containing no great pyramid, and such that no shortest odd hole in $G$ is jewelled.
Let $C$ be a shortest odd hole in $G$, and let $P$ be a path of $G$ with ends $u,v\in V(C)$, such that
$|E(P)|=d_C(u,v)$, and no vertex in $P^*$ is big $C$-major. Let $L_1, L_2$ be the paths of $C$ joining $u,v$,
where $|E(L_1)|=|E(P)|$, and let $C'$ be the shortest odd hole $P\cup L_2$. 
Then every $C'$-major vertex is $C$-major, and vice versa.
\end{thm}
\Proof
Let $P$ have vertices $u\dd p_1\cc p_k\dd v$ in order, and let $C$ have vertices $c_1\cc c_{2n+1}\dd c_1$ in order, where
$u=c_1$ and $v=c_{k+2}$. Define $p_0=u$ and $p_{k+1}=v$.
By \ref{shortcut}, no shortest odd hole has a shortcut.

Suppose that some vertex $w$ of $G$ is $C'$-major and not $C$-major. Since $G$ contains no great pyramid, it follows that
$w$ is big $C'$-major.
Since $w$ is not $C$-major, $w$ has a neighbour
in $P^*$. Choose $h,j\in \{0\ll k+1\}$ minimum and maximum respectively, such that $w$ is adjacent to $p_h, p_j$.

Suppose first that $w\in L_1^*$, and $w=c_i$ say. Then all neighbours of $w$ in $V(C')$ belong to $V(P)$; and there are at least four such neighbours. 
Hence $j\ge h+3$.
Since the path $p_0\dd p_1\cc p_h\dd c_i$ is not a shortcut for $C$, it follows that
$h+1\ge i-1$; and since $c_i\dd p_j\cc p_{k+1}$ is not a shortcut for $C$, $k+2-j\ge k+2-i$. Adding, it follows that
$2\ge j-h$, a contradiction. This proves that $w\notin L_1^*$, and consequently no vertex of $L_1^*$ is $C'$-major.

Hence there is symmetry between $C,C'$, and so if we can prove the first assertion of the theorem, then the ``vice versa''
follows from the symmetry. 

It follows that
$v\notin V(C\cup C')$. Suppose that $j\ge h+2$. Since the path 
$P'$ with vertices $p_0\cc p_h\dd w\dd p_j\cc p_{k+1}$ is not a shortcut for $C$, $j=h+2$; and so $w$ has a neighbour in $L_2^*$, 
since it
has four neighbours in $V(C')$. Hence $P'\cup L_2$ is not a hole, contrary to \ref{reroute}. Thus $j\le h+1$.

If $j=h$, then $w$ has at least three neighbours in $V(L_2)$; and hence exactly three, and they are consecutive, since
$w$ is not $C$-major; and this contradicts \ref{notbig} applied to $C'$. Thus $j=h+1$.

Define $c_{2n+2}=c_1$, and choose $r,t\in \{k+2\ll 2n+2\}$ minimum and maximum respectively, such that $w$ is adjacent to $c_r, c_t$.
Hence $t\le r+2$ since $w$ is not $C$-major; and $t\ge r+1$ since $w$ is big $C'$-major. Also $t\ne r+1$
by \ref{notbig}. Hence $t=r+2$. If $r=k+2$ then $h+1=k+1$, and $C'$ is jewelled, a contradiction; so $r\ge k+3$
and similarly $r+2=t\le 2n+1$. So $w$ is nonadjacent to $u,v$.

Now the path $p_0\dd p_1\cc p_h\dd w\dd c_r$ is not a shortcut for $C$, and so $h+2\ge \min(r-1,2n+2-r)$; but
$h+2<r-1$ since $h+1\le k$, and so $h+2\ge 2n+2-r$. Similarly, by exchanging $u,v$ it follows that $k-h+2\ge r-k$.
Adding, we deduce that $k+1\ge n$. But $L_1$ has length $k+1$ and $L_2$ is longer, so $2(k+1)<2n+1$, and hence $n=k+1$.
Also, $h+2= 2n+2-r$, that is, $h= 2n-r$. The path $p_0\dd p_1\cc p_h\dd w\dd c_r$ has the same length as
the path $c_r\dd c_{r+1}\cc c_{2n+2}$, and less than the path $c_1\dd c_2\cc c_r$; so from \ref{reroute},
the union of the paths $p_0\dd p_1\cc p_h\dd w\dd c_r$ and $c_1\dd c_2\cc c_r$ is a shortest odd hole $C''$ say.
Hence $p_h\dd p_{h+1}\cc p_{k+1}$ is not a shortcut for this hole; and so either one of its internal vertices
is $C''$-major, or $k+1-h\ge d_{C''}(p_h,p_{k+1})$. Now $d_{C''}(p_h,p_{k+1})=\min(h+k+1,r-k)$. Certainly
$k+1-h< h+k+1$ since $h>0$, and $k+1-h< r-k$ since $h= 2n-r$ and $n=k+1$. Thus there exists $i\in \{h+1\ll k\}$ such that
$p_i$ is $C''$-major. But $p_i$ has at most three neighbours in $V(C)$, and so at least one in $V(C')\setminus V(C)$; and so
$i=h+1$. So $p_{h+1}$ has at least two neighbours in the interior of the path $c_1\cc c_r$. By exchanging
$u,v$ it follows similarly that the union of the paths $c_{r+2}\cc c_{2n+1}\dd c_1\cc c_{k+2}$ and $c_{r+2}\dd w\dd p_{h+1}\cc p_{k}\dd c_{k+2}$ is a hole, and so $p_{h+1}$ has no neighbours in $L_1^*$. Hence $p_{h+1}$ has at least two neighbours in the set
$\{c_{k+2}\ll c_r\}$. 

Suppose that $p_{h+1}$ is adjacent to $c_{\ell}$ where $k+3\le \le r-1$.
Now no internal vertex of
the path $p_0\cc p_h\dd p_{h+1}\dd c_{ell}$ is $C$-major; and the length of this path is $h+2$, which is less than
$\min(\ell-1, 2n+2-\ell)$ because $h+1\le k\le \ell-3$ and $h=2n-r< 2n-\ell$. Hence this path is a shortcut for $C$,
a contradiction. Thus $p_{h+1}$ has precisely two neighbours in $\{c_{k+2}\ll c_r\}$, the vertices $c_{k+2}$ and $c_r$.
Consequently $h+1=k$, since $P$ is induced; and $r\le k+4$, since $p_k$ is not $C$-major. 
The path $p_0\cc p_h\dd p_{h+1}\dd c_r$ has the same length as
the path $c_r\cc c_{2n+1}\dd c_1$ and less than the path $c_1\cc c_r$; yet $p_{h+1}$ is adjacent to $c_{k+2}$, 
contrary to \ref{reroute}. This proves \ref{majoreq}.~\bbox

\section{Heavy edges}

We need to use an idea from~\cite{oddholetest}, adapted appropriately. Its proof used other theorems from previous 
papers that assumed there were no pyramids or jewels, so we have to adapt their statements and proofs. Fortunately the changes required 
are very minor, and we think it is unnecessary
to reprint the old proofs in full. 
So we will
just give statements of the theorems we need, which are variants of theorems from~\cite{bergealg} and~\cite{oddholetest},
and sketch how the proofs should be modified. There are two main changes:
\begin{itemize}
\item These old theorems are about $C$-major vertices, and use the fact that when there is no pyramid, all $C$-major vertices are big.
For us, this is not true, since pyramids may be present; but if we just change the statements of the theorems to refer to 
big $C$-major vertices then all is well.
\item These old theorems assume that $G$ contains no jewel, but in every case, all that the proof needs is that no shortest odd hole 
is jewelled, and this we can assume.
\end{itemize}
Let us see these old theorems in detail. 

Theorem 7.6 of~\cite{bergealg} assumes that $G$ has no pyramid or jewel (and its proof uses theorem 7.5 of~\cite{bergealg}, which
we have to abandon), 
but the changes of the bullets above repair its proof.
Since ``normal'' subsets 
are nonempty, these modifications allow us to prove:
\begin{thm}\label{7.6mod}
Let $G$ be a graph in which no shortest odd hole is jewelled. Let $C$ be a shortest odd hole in $G$, and let $X$ be a stable set
of big $C$-major vertices. Then there is a vertex $v\in V(C)$ adjacent to every vertex in $X$.
\end{thm}

We need nothing more from~\cite{bergealg}, but we need some results from~\cite{oddholetest}. Theorem 3.3 of~\cite{oddholetest}
assumes that $G$ has no pyramid or jewel, and its proof uses theorem  7.6 of~\cite{bergealg}.
But again, the changes of the bullets 
fix both these problems: use \ref{7.6mod} above in place of theorem 7.6 of~\cite{bergealg},
and assume that no shortest odd hole is jewelled.
We obtain:
\begin{thm}\label{3.3mod}
Let $G$ be a graph in which no shortest odd hole is jewelled. Let $C$ be a shortest odd hole in $G$, and let $x,y$ be
nonadjacent big $C$-major vertices.
Then every induced path between $x, y$ with interior in $V(C)$ has
even length.
\end{thm}

Theorem 3.4 of~\cite{oddholetest} assumes that $G$ has no pyramid or jewel, but the changes of the bullets repair that; 
and its proof uses theorem 3.3 of~\cite{oddholetest}.
but we can replace it by \ref{3.3mod} above.
We obtain:
\begin{thm}\label{3.4mod}
Let $G$ be a graph with no hole of length five, and in which no shortest odd hole is jewelled, and let $C$ be a shortest odd hole in $G$.
Let $X$ be a set of big $C$-major vertices, and let $x_0 \in X$ be nonadjacent to all other members of $X$.
Then there is an edge $uv$ of $C$ such that every member of $X$ is adjacent to one of $u,v$.
\end{thm}

\section{The proof of \ref{nogreatpyr}.}
In this section we prove \ref{nogreatpyr}. It is almost identical with the algorithm of section four of~\cite{oddholetest}.
First, we need a version of theorem 4.2 of~\cite{bergealg}, the following:

\begin{thm}\label{testclean}
There is an algorithm with the following specifications:
\begin{description}
\item [Input:] A graph $G$.
\item [Output:] Either an odd hole of $G$, or a statement of failure. If $G$ contains no jewelled shortest odd hole and no great pyramid,
and some shortest odd hole in $G$ is clean,
the output is a shortest odd hole of $G$.
\item [Running time:] $O(|V(G)|^4)$.
\end{description}
\end{thm}
\Proof Here is an algorithm. For every pair of vertices $u,v$, find a shortest path $P(u,v)$ between them,
if one exists. For every triple $u,v,w$, test whether the three paths $P(u,v),P(v,w),P(w,u)$
all exist, and if so whether their union is an odd hole. If we find such a hole, record it.
When all triples have been examined, if no hole has been recorded, report failure, and otherwise output the recorded hole
of smallest length.

To see that this works correctly, we may assume that contains no jewelled shortest odd hole and no great pyramid,
and some shortest odd hole $C$ in $G$ is clean.
Choose vertices $u,v,w \in V(C)$, roughly equally spaced
in $C$; more precisely, such that every component of $C\setminus\{u,v,w\}$ contains at most
$n-1$ vertices, where $C$ has length $2n+1$.
Since there is a path joining $u,v$, the algorithm will find
a shortest such path $P(u,v)$. We claim that $C$ can be chosen containing $P(u,v)$. By \ref{shortcut}, there is no
shortcut for $C$, since $G$ contains no great pyramid. Let $L_1$ be the path of $C$
joining $u,v$, not passing through $w$.
Then $L_1$ has length $\le n$, from the choice of $u,v,w$, and so since $P(u,v)$ is not a shortcut for $C$, and none of
its vertices are $C$-major since $C$ is clean, it follows that 
$L_1,P(u,v)$ have the same length. Let $L_2$ be the
second path of $C$ between $u,v$ in $C$. The union of $L_2,P(u,v)$ is a clean shortest odd hole, by \ref{reroute} and \ref{majoreq},
and so we may choose $C$ containing $P(u,v)$. By repeating this for the other two pairs from $u,v,w$,
we see that $C$ can be chosen to include all of $P(u,v),P(v,w),P(w,u)$
simultaneously. So the union of the three paths joining $u,v,w$ chosen by the algorithm is a shortest odd hole,
and therefore in this case the algorithm correctly records a shortest odd hole, and therefore will output one.

The running time of the algorithm as described is
$O(|V(G)|^5)$, because after selecting $u,v,w$ and the three paths, it takes quadratic time to check
whether the three paths make a hole. Here is a (sketch of) how to get the running time down to $O(|V(G)|^4)$, although
it makes no difference to the running time of our main algorithm. For each pair of distinct vertices $u,v$, mark the vertices that 
belong to $P(u,v)$ or have a neighbour in its interior; then for all $w$, we can compute in linear time whether $P(u,v)\cup P(v,w)$
is an induced path,
and whether $P(w,u)\cup P(u,v)$ is an induced path. Then with all this information (which takes time $O(|V(G)|^4)$ to compute)
we check whether there is a triple $u,v,w$ of distinct vertices such that each pair of the paths $P(u,v), P(v,w), P(w,u)$ makes an induced path; and
checking each triple now takes constant time.
This proves \ref{testclean}.\bbox

Let us say a shortest odd hole $C$ is {\em heavy-cleanable} if there is an edge $uv$ of $C$ such that
every big $C$-major vertex is adjacent to one of $u,v$. We deduce:

\begin{thm}\label{testcleanable}
There is an algorithm with the following specifications:
\begin{description}
\item [Input:] A graph $G$.
\item [Output:] Either an odd hole of $G$, or a statement of failure. If $G$ contains no 5-hole, no jewelled shortest odd hole
and no great pyramid, and contains a heavy-cleanable shortest odd hole, the output is a shortest odd hole.
\item [Running time:] $O(|G|^8)$.
\end{description}
\end{thm}
\Proof
List all the four-vertex induced paths $c_1\dd c_2\dd c_3\dd c_4$ of $G$. For each one, let $X$ be the set of all vertices of $G$
different from $c_1,\ldots, c_4$ and adjacent to one of $c_2,c_3$. We run \ref{testclean} on $G\setminus X$ and 
record any hole that it outputs.
If after examining all 4-tuples, no hole is recorded, report failure, and otherwise output the recorded hole of smallest length.

To see correctness, we may assume that $G$ has no 5-hole, no jewelled shortest odd hole
and no great pyramid, and contains a heavy-cleanable shortest odd hole $C$. Thus
$C$ is clean in $G\setminus X$ for some $X$ that we test;
and when we do so,
\ref{testclean} outputs a shortest odd hole of $G$, that we record. Consequently in this case the algorithm outputs a shortest
odd hole of $G$.
This proves \ref{testcleanable}.~\bbox

\begin{thm}\label{noheavyclean}
There is an algorithm with the following specifications:
\begin{description}
\item [Input:] A graph $G$.
\item [Output:] Outputs either an odd hole of $G$, or a statement of failure. If $G$ contains no 5-hole, no jewelled shortest 
odd hole, no great pyramid, and no heavy-cleanable shortest odd hole, the output will
be a shortest odd hole of $G$.
\item [Running time:] $O(|G|^9)$.
\end{description}
\end{thm}
\Proof  This is the algorithm described in section 4 of~\cite{oddholetest}, using \ref{shortcut} and \ref{reroute}
in place of theorem 4.1 of~\cite{oddholetest}, and using \ref{3.4mod} in place of theorem 3.4 of~\cite{oddholetest}.~\bbox

Now we are ready to prove \ref{nogreatpyr}, which we restate:
\begin{thm}\label{nogreatpyragain}
There is an algorithm with the following specifications:
\begin{description}
\item [Input:] A graph $G$.
\item [Output:] Outputs either an odd hole of $G$, or a statement of failure. If $G$ contains no 5-hole, no jewelled shortest odd hole, 
and no
great pyramid, and $G$ contains an odd hole, the output will
be a shortest odd hole of $G$.
\item [Running time:] $O(|G|^9)$.
\end{description}
\end{thm}
We are given an input graph $G$. We apply the algorithms of \ref{testcleanable} and \ref{noheavyclean} to $G$, and record the holes
that they output. If no hole is recorded, report failure, and otherwise output the shortest recorded hole.

To see correctness, we may assume that $G$ contains no 5-hole, no jewelled shortest odd hole, 
and no
great pyramid. If $G$ has a heavy-cleanable shortest odd hole, \ref{testcleanable} will output a shortest odd hole. 
If $G$ has no odd hole, both algorithms will report failure, and the output is correct. Otherwise $G$ has an odd hole
and has no heavy-cleanable shortest odd hole; and then \ref{noheavyclean} will output a shortest odd hole. 
This proves \ref{nogreatpyragain} and hence \ref{nogreatpyr}.~\bbox

\section{Cleaning a great pyramid}

It remains to prove \ref{findgreatpyr}, but that requires several lemmas. 
Let $H$ be a great pyramid in $G$, with apex $a$ and base $\{b_1,b_2,b_3\}$, and 
constituent paths $P_1, P_2, P_3$ where $P_i$ is between $a,b_i$, and $P_3$ is shorter than $P_1, P_2$. Thus $G[V(P_1\cup P_2)]$
is a shortest odd hole, $C$ say. 
We call the length of $P_3$ the {\em height} of the pyramid. 

The idea of the algorithm for \ref{findgreatpyr} is, we look for a great pyramid with minimum height, height $r$ say; then 
because of \ref{shortcut}, we know the useful fact that no shortest odd hole of $G$ has a shortcut of length at most $r$. We will guess 
a few important
vertices of the great pyramid, and then try to fill in the paths between them by picking shortest paths in appropriate subgraphs, 
using the ``useful fact''. But shortcuts by definition contain no big $C$-major vertices,
so to use the useful fact, we need to be sure there are no big $C$-major vertices in our paths, which is tricky because we do 
not know $C$. Ideally we would first clean
to get rid of all big $C$-major vertices, but it turns out that a partial cleaning via \ref{3.4mod} is enough, because of a convenient property of 
big $C$-major vertices that we prove in \ref{pyrmajor} below.

For brevity, let us say a vertex $v$ is {\em major} for the great pyramid $H$ if it is big $C$-major.
Let $v$ be major for $H$, and let $\{i,j,k\}= \{1,2,3\}$. 
We say that $v$ has {\em type} $(P_i, P_j)$ if
\begin{itemize}
\item $v$ has at least three neighbours in $V(P_i)\setminus \{a\}$;
\item $v$ has exactly two neighbours in $V(P_j)$ and they are adjacent; and
\item $v$ has no neighbours in $V(P_k)\setminus \{a\}$.
\end{itemize}
\begin{thm}\label{pyrmajor}
Let $H$ be a great pyramid, with notation as above, and let $v$ be major for $H$. Then
either $v$ has at least two neighbours in $\{b_1,b_2,b_3\}$, or 
$v$ has type $(P_i, P_j)$ where $(i,j)$ is one of the pairs $(1,2), (2,1), (1,3), (2,3)$.
\end{thm}
\Proof We may assume that $v$ has at most one neighbour in $\{b_1,b_2,b_3\}$.
Let $\ell_i=|E(P_i)|$ for $i = 1,2,3$.
\\
\\
(1) {\em There exists $k\in \{1,2,3\}$
such that $v$ has no neighbour in $V(P_k)\setminus \{a\}$.}
\\
\\
Suppose that $v$ has a neighbour in $V(P_i)\setminus \{a\}$ for all $i\in \{1,2,3\}$.
For $i = 1,2,3$, choose a minimal path $P_i'$ between $v$ and $b_i$ with interior in $P_i^*$. 
Let $|E(P_i')|=\ell_i'$ for $i = 1,2,3$.
Then at least two of $\ell_1', \ell_2',\ell_3'\ge 2$, so these three paths are the constituent paths of a pyramid. 
Hence some two of them induced an odd hole
$C'$, say $P_i', P_j'$. Since $C'$ is not shorter than $C$, it follows that $\ell_i'+\ell_j'\ge \ell_1+\ell_2$.
But $\ell_i'\le \ell_i$ for $i = 1,2,3$, and hence $\ell_i+\ell_j\ge \ell_1+\ell_2$. Since $\ell_3<\ell_1, \ell_2$,
it follows that $\{i,j\}=\{1,2\}$; and $P_1'$ has the same length as $P_1$, so $v$ has a unique neighbour in $V(P_1)\setminus \{a\}$,
the neighbour of $a$ in $P_1$. The same holds for $P_2$; but then $v$ has at most three neighbours in $V(C)$, and they all belong to 
a three-vertex subpath of $C$, contradicting that $v$ is $C$-major. This proves (1).
\\
\\
(2) {\em $v$ has neighbours in exactly two of the sets $V(P_i)\setminus \{a\}\;(i = 1,2,3)$.}
\\
\\
Suppose that there is at most one such $i$; so all neighbours of $v$ in $V(H)$ belong to $V(P_i)$. 
So $i\in \{1,2\}$ (because $v$ is $C$-major), and we may assume that $i = 1$. Thus $v$ has at least three neighbours in $V(P_1)$,
Let $P_1'$ be the induced path between $a, b_1$
with interior in $P_1^*\cup \{v\}$ that contains $v$, with length $\ell_1'$ say. So $\ell_1'<\ell_1$ because the neighbours
of $v$ in $V(P_1)$ do not all lie in a three-vertex subpath; and $\ell_1'\ge 2$, so the three paths $P_1', P_2,P_3$
define a pyramid, and hence some two of $P_1', P_2,P_3$ induce an odd hole $C'$. But every two of $\ell_1', \ell_2, \ell_3$
sum to less than $\ell_1+\ell_2$, since $\ell_1'<\ell_1$ and $\ell_3<\ell_1, \ell_2$, a contradiction. This proves (2).

\bigskip

Since $v$ is big $C$-major, it has two nonadjacent neighbours in one of $V(P_1), V(P_2)$, by \ref{notbig}; and by
exchanging $P_1,P_2$ if necessary, we may assume that $v$ has two nonadjacent neighbours in $V(P_1)$.
By (1), (2) there exists a unique $j\in \{2,3\}$ such that $v$ has a neighbour in $V(P_j)\setminus \{a\}$.
Let $R_1$ be the induced path between $v, a$ with interior in $P_1^*$, and let $S_1$ be the induced path between $v, b_1$
with interior in $P_1^*$. Define $R_j, S_j$ similarly. Thus $V(R_1)\cap V(S_1)=\{v\}$, and $R_1\cup S_1$
is an induced path between $a, b_1$, but this need not be true for $R_j, S_j$ since $v$ might not have two nonadjacent
neighbours in $V(R_j)$. Let $\{j,k\}=\{2,3\}$.
\\
\\
(3) {\em We may assume that $v$ has at least two neighbours in $V(P_j)$.}
\\
\\
For suppose that $v$ has a unique neighbour $u\in V(P_j)$. Since $v$ has a neighbour in $V(P_j)\setminus \{a\}$, it follows
that $u\ne a$, and $v,a$ are nonadjacent.

Suppose first that $j=2$. 
There is a pyramid with apex $u$ and
constituent paths
$u\dd v\dd S_1\dd b_1$, $(R_2\setminus v)\cup P_3$, and $S_2\setminus v$. Some two of these three paths induce an odd hole,
and so the sum of their lengths is at least $\ell_1+\ell_2$. The last two sum to $\ell_2+\ell_3<\ell_1+\ell_2$;
and the first and third sum to $|E(S_1)|+|E(S_2)|<\ell_1+\ell_2$ since $|E(S_1)|<\ell_1-3$ and $|E(S_2)|\le \ell_2$.
Thus the sum of the lengths of $u\dd v\dd S_1\dd b_1$ and  $(R_2\setminus v)\cup P_3$ is at least $\ell_1+\ell_2$; and so
$$|E(S_1)|+|E(R_2)|+\ell_3\ge \ell_1+\ell_2.$$
In particular, $|E(R_2)|$ has length at least three, since $|E(S_1)|\le \ell_1-2$ and $|E(P_3)|\le \ell_2-1$.
So $u$ is nonadjacent to $a$.
On the other hand, there is a pyramid with apex $v$ and constituent paths $S_1, R_1\cup P_3$, and $S_2$; so some two of the lengths
of these sum to at most $\ell_1+\ell_2$. But 
the first two sum to at most $\ell_1+\ell_3<\ell_1+\ell_2$, and the first and third
sum to at most $(\ell_1-2)$, since $S_1$ has length at most $\ell_1-2$; so the lengths of $R_1\cup P_3$ and $S_2$
sum to at least $\ell_1+\ell_2$. Thus 
$$|E(R_1)|+\ell_3+|E(S_2)|\ge \ell_1+\ell_2.$$
Summing with the previous displayed inequality, we deduce
$$|E(S_1)|+|E(R_2)|+ |E(R_1)|+2\ell_3+|E(S_2)|\ge 2(\ell_1+\ell_2).$$
But $|E(R_1)|+|E(S_1)|\le \ell_1$, and $|E(R_2)|+|E(S_2)|= \ell_2+2$, and so
$$2+ 2\ell_3\ge \ell_1+\ell_2.$$
It follows that $\ell_3+1=\ell_1=\ell_2$, since $\ell_3<\ell_1, \ell_2$; and equality holds throughout. In particular,
$|E(R_1)|+|E(S_1)|=\ell_1$, and so $v$ has only three neighbours in $V(P_1)$, and they are consecutive. But this contradicts
\ref{notbig}.

This proves that $j=3$. Here are three holes, all shorter than $C$ and hence even:
$$R_1\cup R_3;\; v\dd S_1\dd b_1\dd b_3\dd S_3\dd v;\; a\dd P_2\dd b_2\dd b_3\dd P_3\dd a.$$
(To see that the second is a hole we use that $v$ has at most one neighbour
in $\{b_1,b_2,b_3\}$.) Let us add a fourth cycle to this list, the odd cycle $b_1\dd b_2\dd b_3\dd b_1$.
The symmetric difference of the edge sets of these four cycles is the hole 
$$a\dd P_2\dd b_2\dd b_1\dd S_1\dd v\dd R_1\dd a,$$
and so the latter is odd, and yet it is shorter than $C$ (because $v$ is $C$-major), a contradiction. This proves (3).
\\
\\
(4) {\em $v$ does not have two nonadjacent neighbours in $V(P_j)$.}
\\
\\
Suppose it does. Again, there are two cases, depending whether $j=2$ or $3$. Suppose first that $j=2$. There is a pyramid with apex $v$
and constituent paths $S_1, S_2, R_1\cup P_3$, and by the usual argument we deduce that 
$$|E(R_1)+\ell_3+|E(S_2)|\ge \ell_1+\ell_2.$$
From the symmetry between $P_1, P_2$, it also follows that
$$|E(R_2)+\ell_3+|E(S_1)|\ge \ell_1+\ell_2.$$
Adding, we deduce that
$$|E(R_1)|+|E(S_1)|+|E(R_2)|+|E(S_2)|+2\ell_3\ge 2(\ell_1+\ell_2).$$
Since $|E(R_i)|+|E(S_i)|\le \ell_i$ for $i = 1,2$, it follows that
$2\ell_3\ge \ell_1+\ell_2$, a contradiction.

Thus $j=3$. The pyramid with apex $v$ and constituent paths $S_1,R_1\cup P_2, S_3$ tells us that
$$|E(R_1)|+\ell_2+|E(S_3)|\ge \ell_1+\ell_2.$$
But the pyramid with apex $v$ and constituent paths
$S_1,R_3\cup P_2, S_3$ tells us that 
$$ |E(S_1)|+|E(R_3)|+\ell_2\ge \ell_1+\ell_2.$$
Adding, and using that $|E(R_i)|+|E(S_i)|\le \ell_i$ for $i = 1,3$, we deduce that
$\ell_1+2\ell_2+\ell_3\ge 2(\ell_1+\ell_2)$, which simplifies to $\ell_3\ge \ell_1$, a contradiction.  This proves (4).

\bigskip
From (3) and (4), $v$ has exactly two neighbours in $V(P_j)$, and they are adjacent. Thus $v$ has type $(P_1, P_j)$. This
proves \ref{pyrmajor}.~\bbox

\section{Jumps off a great pyramid}

Let us use the same notation as in the last section: thus $H$ is a great pyramid, with apex $a$, and constituent paths $P_1,P_2,P_3$,
where $P_i$ has ends $a, b_i$ for $i = 1,2,3$, and $P_3$ is shorter than $P_1,P_2$. let $C$ be the shortest 
odd hole $G[V(P_1\cup P_2)]$, and let $P_i$ have length $\ell_i$ for $i = 1,2,3$. 
We remark that $\ell_1,\ell_2$ have the same parity (since $C$ has odd length), and $\ell_3$ has the other parity (since otherwise
the hole $G[V(P_1\cup P_3)]$ would have odd length and be shorter than $C$).
Let us say $H$ is {\em optimal} if there is no great pyramid
in $G$ with smaller height. 
In this section we develop some results about short paths in $G$ that join vertices in $V(H)$, when $H$ is optimal.
We begin with:
\begin{thm}\label{optshortcut}
Let $H$ be an optimal great pyramid in $G$, with height $\ell_3$. Then no shortest odd hole in $G$ has a shortcut of length at most 
$\ell_3$.
\end{thm}
\Proof
Suppose there is such a shortcut, $P$ say; then by \ref{shortcut}, there is a great pyramid such that all vertices of 
its shortest
constituent path (except the apex) belong to $P^*$, and so there are at most $\ell_3-1$ of them. Hence the great pyramid
has height at most $\ell_3-1$, contradicting that $H$ is optimal. This proves \ref{optshortcut}.~\bbox

\begin{thm}\label{chooseP3}
Let $H$ be an optimal great pyramid in $G$. In the notation above let $X$ be the set of all big $C$-major vertices,
together with all vertices adjacent or equal to $b_1$ or to $b_2$. 
Let $P_3'$ be a shortest path between $a,b_3$ with interior in $V(G)\setminus X$. 
Then $P_3',P_3$ have the same length, and the subgraph
induced on $V(P_1\cup P_2\cup P_3')$ is an optimal great pyramid.
\end{thm}
\Proof
Since no vertex in $P_3^*$ is $C$-major, it follows that $P_3^*$ contains no vertex in $X$, and so $P_3'$ has length at most that of $P_3$.
We claim that $P_1,P_2,P_3'$ are the constituent paths of a pyramid. If so, then it is a great pyramid, and so $P_3,P_3'$ have
the same length since $H$ is optimal. So we only need to show the claim that $P_1,P_2,P_3'$ are the constituent paths of a pyramid.
Suppose not; then some vertex of $P_3'^*$ belongs to or has a neighbour in one of $V(P_1)\setminus \{a\}$, $V(P_2)\setminus \{a\}$.
Choose a minimal subpath $Q$ of $P_3'$ between $b_3$ and some vertex $q$ with a neighbour in one of 
$V(P_1)\setminus \{a\}$, $V(P_2)\setminus \{a\}$. From the symmetry we may assume that $q$ has a neighbour in 
$V(P_1)\setminus \{a\}$. 
\\
\\
(1) {\em $q$ has no neighbour in $V(P_2)\setminus \{a\}$.}
\\
\\
Suppose it does. Then since $q$ is nonadjacent to $b_1,b_3$, there is a pyramid with apex $q$
and constituent paths the path $Q$, and for $i = 1,2$ an induced path between $q,b_i$ with interior in $P_i^*$. 
Since $Q$ has length less than 
$\ell_3$, and the other two constituent paths have lengths at most $\ell_1$, $\ell_2$ respectively, this contradicts that
$H$ is an optimal great pyramid. This proves (1).

\bigskip

No vertex of $Q$ is big $C$-major. Let $p$ be the neighbours of $q$ in $P_1$ that is closest in $P_1$ to $b_1$. Let $R_1,S_1$
be the two subpaths of $P_1$, between $p,a$ and between $p, b_1$ respectively. 
\\
\\
(2) {\em Either $|E(Q)|+1=\ell_3$, or $|E(Q)|\ge |E(S_1)|-1$.}
\\
\\
Certainly $|E(Q)|\le \ell_3-1$, since $Q$ is a proper subpath of $P_3'$; suppose that $|E(Q)|< \ell_3-1$. Then the path
$b_2\dd b_3\dd Q\dd q\dd p$ has length at most $\ell_3$, and so is not a shortcut for $C$. Hence its length is at least
$d_C(b_2,p)$. Thus either $|E(Q)|+2\ge |E(S_1)|+1$, or $|E(Q)|+2\ge |E(R_1)|+\ell_2$. Suppose the second holds.
Since $|E(Q)|\le \ell_3-1\le \ell_2-2$, it follows that $E(R_1)|=0$, a contradiction. So the first holds. This proves (2).

\bigskip

Let $Q'$ be the subpath of $P_3'$ between $q, a$.
\\
\\
(3) {\em $|E(Q')|\ge |E(R_1)|-1$.}
\\
\\
Suppose not. The path $p\dd q\dd Q'\dd a$ is has length at most $\ell_3$, since $|E(Q')|\le |E(P_3')|-1\le \ell_3-1$, and
so it is not a shortcut for $C$. Hence its length is at least $d_C(p,a)$. But $E(Q')|+1\le \ell_3<\ell_2$, and so
Thus $|E(Q')|+1\ge |E(R_1)|$. This proves (3).
\\
\\
(4) {\em $|E(Q)|+1<\ell_3$.}
\\
\\
Suppose that $|E(Q)|+1=\ell_3$. Then $q$ is adjacent to $a$, since $Q$ is a subpath of $P_3'$; and so $|E(Q')|=1$. From (3)
(or since $q$ is not $C$-major) it follows that $|E(R_1)|\le 2$. Suppose that $|E(R_1)|=2$. Then the hole 
$$a\dd q\dd p\dd S_1\dd b_1\dd b_2\dd P_2\dd a$$ 
is a shortest odd hole, and since the path $b_2\dd b_1\dd Q\dd q$ has length at most $\ell_3$, it is not a shortcut over this hole.
Since $b_2\dd b_1\dd Q\dd q$ has length at most $\ell_3<\ell_2$, it follows that $|E(Q)|\ge |E(S_1)|+1=\ell_1-1$,
which is impossible since $\ell_3<\ell_1$. Thus $|E(R_1)|=1$. But then the hole 
$$b_3\dd Q\dd q\dd p\dd S_1\dd b_1\dd b_2$$
has length $\ell_1+\ell_3$, and since this is less than the length of $C$, and $\ell_1,\ell_3$ have opposite parity, this is
a contradiction.
This proves (4).
\\
\\
(5) {\em 
$|E(Q')|=|E(R_1)|-1$.}
\\
\\
From (2) and (4) we deduce that $|E(Q)|\ge |E(S_1)|-1$, and adding the inequality of (3), it follows that
$|E(P_3')|\ge \ell_1-2$. 
But $|E(P_3')|\le \ell_3$, and so $\ell_3\ge \ell_1-2$. Since
$\ell_3<\ell_1$ and $\ell_1, \ell_3$ have opposite parity, it follows that $\ell_3=\ell_1-1$. 
We have 
$$(|E(P_3')|-\ell_3) +(\ell_3-\ell_1)+2 = (|E(Q)|-(|E(S_1)|-1)) + (|E(Q')|-(|E(R_1)|-1)).$$
Since $\ell_3=\ell_1-1$, we deduce that
$$(\ell_3-|E(P_3')|)+ (|E(Q)|-(|E(S_1)|-1)) + (|E(Q')|-(|E(R_1)|-1))=1.$$
In particular, $|E(Q)|-(|E(S_1)|-1)\le 1$. 
But $Q,S_1$ have lengths of the same parity, since the subgraph induced on $V(Q\cup S_1)$ is a hole of length less than $C$;
so $|E(Q)|=|E(S_1)|$, and hence
$$(\ell_3-|E(P_3')|)+ (|E(Q')|-(|E(R_1)|-1))=0.$$
Thus $P_3'$ has length $\ell_3$, and $|E(Q')|=|E(R_1)|-1$. This proves (5).

\bigskip

By (5), the path $p\dd q\dd Q'\dd a$ has the same length as $R_1$. 
Since the length of the path $p\dd q\dd Q'\dd a$ is at most $\ell_3$, no proper subset of its interior is a shortcut
for any shortest odd hole; so by \ref{reroute}, either
$$a\dd P_2\dd b_2\dd b_1\dd S_1\dd p\dd q\dd Q'\dd a$$
is a shortest odd hole $C'$ say, or $G[V(C\cup Q')]$ is a great pyramid with heart $V(Q')\setminus \{a\}$. In the first case, 
the path $b_2\dd b_1\dd Q\dd q$ is a shortcut for $C'$, since none of its internal vertices are big $C'$-major; and yet 
its length is $|E(Q)|+1\le \ell_3$, a contradiction. So the second holds, and $G[V(C\cup Q')]$ is a great pyramid $H'$ with heart 
$V(Q')\setminus \{a\}$. But $|V(Q')\setminus \{a\}|<|V(P_3)\setminus \{a\}|$, so the height of $H'$ is less than the height
of $H$, contrary to the optimality of $H$. This proves \ref{chooseP3}.~\bbox

\begin{thm}\label{chooseR2}
Let $H$ be an optimal great pyramid in $G$. In the notation as before let $X$ be the set of all big $C$-major vertices,
together with all vertices adjacent or equal to a vertex in $\{b_1\}\cup (V(P_3)\setminus \{a\})$.
Let $c_2\in V(P_2)$ such that the subpath $R_2$ of $P_2$ between $c_2,a$ has length at most $\min(\ell_3, (1+\ell_2)/2)$.
Let $R_2'$ be a shortest path between $a,c_2$ with interior in $V(G)\setminus X$. Then $R_2',R_2$ have the same length, and the subgraph
induced on $V(P_1\cup P_3)\cup (V(P_2)\setminus V(R_2))\cup V(R_2')$ is an optimal great pyramid.
\end{thm}
\Proof
Since no vertex of $R_2^*$ belongs to $X$, it follows that the length of $R_2'$ is at most that of $R_2$. 
\\
\\
(1) {\em We may assume (for a contradiction) that the interior of $R_2'$ is not anticomplete to $P_1^*$.}
\\
\\
Suppose it is;  then it is also anticomplete to $(V(P_1)\cup V(P_2))\setminus \{a\}$. Choose 
an induced path $P_2'$ between $c_2,a$ with interior in $(V(P_2)\setminus V(R_2))\cup V(R_2')$, we see that $P_1,P_2',P_3$ are the
constituent paths of a pyramid $H'$, and so some two of these three paths have sum of lengths at least $\ell_1+\ell_2$. But $P_2'$
has length at most $\ell_2$, and $P_3$ has length $\ell_3<\ell_1,\ell_2$, so the sum of the lengths of $P_1,P_2'$
is at least $\ell_1+\ell_2$. Consequently $P_2'$ has length at least $\ell_2$, and so exactly $\ell_2$. Hence $R_2'$
has the same length as $R_2$, and $H'$ is a great pyramid, and the theorem holds.
This proves (1).

\bigskip

From (1) we may choose a minimal subpath $Q$ of $P_2'$, with ends $c_2,q$ say, such that $q$ has a neighbour in $P_1^*$.
Choose $p\in P_1^*$ adjacent to $q$ such that the subpath of $P_1$ between $p, b_1$ is minimal, and let $S_1$ be this
subpath. Let $R_1$ be the subpath of $P_1$ between $p$ and $a$. Now no vertex in $V(Q)$ is big $C$-major; and
$$d_C(c_2,p)=\min(|E(R_2)|+|E(R_1)|,|E(S_1)|+1+\ell_2-|E(R_2)|).$$
Since $|E(R_2)|+|E(R_1)|> |E(Q)|+1$ (because $R_2$ has length more than $Q$), and $|E(S_1)|+1+\ell_2-|E(R_2)|>|E(Q)|+1$
(because $|E(S_1)|\ge 1$ and $\ell_2-|E(R_2)|>    |E(R_2)|-2\ge |E(Q)|-1$), it follows that
the path $c_2\dd Q\dd q\dd p$ is a shortcut for $C$. But every shortcut has length
at least $\ell_3+1$, from \ref{shortcut} and the optimality of $H$, and yet $Q$ has length less than $\ell_3$,
a contradiction. This proves \ref{chooseR2}.~\bbox

\begin{thm}\label{chooseS2}
Let $H$ be an optimal great pyramid in $G$. In the notation as before let $X$ be the set of all big $C$-major vertices,
together with all vertices adjacent or equal to a vertex in $\{b_1\}\cup (V(P_3)\setminus \{a\})$.
Let $c_2\in V(P_2)$ such that the subpath $S_2$ of $P_2$ between $c_2,b_2$ has length at most $\min(\ell_3, \ell_2/2)$.
Let $S_2'$ be a shortest path between $c_2, b_2$ with interior in $V(G)\setminus X$. Then $S_2',S_2$ have the same length, 
and the subgraph
induced on $V(P_1\cup P_3)\cup (V(P_2)\setminus V(S_2))\cup V(S_2')$ is an optimal great pyramid.
\end{thm}
\Proof As in the proof of \ref{chooseR2}, we may assume that $S_2'^*$ is not anticomplete to $P_1^*$, and so 
we can choose a minimal subpath $Q$ of $P_2'$, with ends $c_2,q$ say, such that $q$ has a neighbour in $P_1^*$.
Choose $p\in P_1^*$ adjacent to $q$ such that the subpath of $P_1$ between $p, b_1$ is minimal, and let $S_1$ be this
subpath. Let $R_1$ be the subpath of $P_1$ between $p$ and $a$. Now no vertex in $V(Q)$ is big $C$-major; and
$$d_C(c_2,p)=\min(|E(R_1)|+\ell_2-|E(S_2)|,|E(S_1)|+1+|E(S_2)|).$$
Since $|E(R_1)|+\ell_2-|E(S_2)|> |E(Q)|+1$ (because $|E(R_1)|\ge 1$ and $\ell_2-|E(S_2)|\ge |E(S_2)| >|E(Q)|$),
and $|E(S_1)|+1+|E(S_2)|>|E(Q)|+1$ (because $|E(S_1)|\ge 1$ and $|E(S_2)|>|E(Q)|$ )
it follows that
the path $c_2\dd Q\dd q\dd p$ is a shortcut for $C$. But every shortcut has length
at least $\ell_3+1$, from \ref{shortcut} and the optimality of $H$, and yet $Q$ has length less than $\ell_3$,
a contradiction. This proves \ref{chooseR2}.~\bbox

\begin{thm}\label{chooseC2D2}
Let $H$ be an optimal great pyramid in $G$. In the notation as before let $X$ be the set of all big $C$-major vertices,
together with all vertices adjacent or equal to a vertex in $\{b_1\}\cup V(P_3)$.
Suppose that $\ell_2\ge 2\ell_3$; and let $R_2, S_2$ be subpaths of $P_2$ both of length $\ell_3$, with one end $a$ and $b_2$
respectively. Let their other ends be $c_2,d_2$ respectively. Let $m_2$ be the vertex of $P_2$ such that the subpath of $P_2$
between $m_2,a$ has length $\lceil \ell_2/2\rceil$. Let $C_2, D_2$ be the subpaths of $P_2$ between $m_2,c_2$ and between $m_2,d_2$
respectively. 
Let $C_2'$ be a shortest path between $m_2,c_2$ with interior in $V(G)\setminus X$, and let $D_2'$ be a shortest path
between $m_2,c_2$ with interior in $V(G)\setminus X$. Then $C_2', C_2$ have the same length, and $D_2', D_2$
have the same length, and 
and the subgraph
induced on $V(P_1\cup P_3)\cup (V(P_2)\setminus V(C_2\cup D_2))\cup V(C_2'\cup D_2')$ is an optimal great pyramid.
\end{thm}
\Proof
Suppose not; then as before, some vertex of one of $C_2',D_2'$ has a neighbour in $P_1^*$, and 
we can choose a minimal subpath $Q$ of one of $C_2', D_2'$, with ends $m_2,q$ say, such that $q$ has a neighbour in 
$P_1^*$. Thus $Q$ has length at most $\lceil \ell_2/2\rceil -\ell_3 -1$. 
Choose an induced path $P_2'$ between $q, b_2$ with interior in $V(Q\cup S_2)$. Thus $P_2'$ has length at most
$$|E(Q)|+\lfloor \ell_2/2 \rfloor\le \ell_2 -\ell_3 -1\le  \ell_2 -\ell_3-1.$$
Now there are three cases, depending whether
$q$ has one neighbour, two nonadjacent neighbours, or just two adjacent neighbours in $V(P_2)$.
\\
\\
(1) {\em $q$ does not have a unique neighbour in $V(P_1)$.}
\\
\\
Suppose it does, $p$ say. Let $R_1$ be the subpath of $P_1$ between $p$ and $a$, and let $S_1$ be the subpath of $P_1$ between $p, b_1$.
There is a pyramid with apex $p$ and constituent paths $S_1, R_1\cup P_3$ and $p\dd q\dd P_2'\dd b_2$, so some two of these paths
have sum of lengths at most $\ell_1+\ell_2$. The first two lengths sum to $\ell_1+\ell_3<\ell_1+\ell_2$; and the first and third
sum to at most 
$$(\ell_1-1) + (\ell_2 -\ell_3) <\ell_1+\ell_2;$$
and the last two sum to at most 
$$(|E(R_1|+\ell_3)+ (\ell_2 -\ell_3)< \ell_1+\ell_2,$$
a contradiction. This proves (1).
\\
\\
(2) {\em $q$ does not have two nonadjacent neighbours in $V(P_1)$.}
\\
\\
Suppose it does. Since $q$ is not $C$-major, there is a three-vertex subpath of $P_1$, with ends $r_1,s_1$, where $a,r_1,s_1,b_1$
are in order in $P_1$, such that $q$ is adjacent to $r_1,s_1$, and possible the vertex of $P_1$ between them,
and has no other neighbours in $V(P_1)$. Thus $s_1\ne b_1$ since $q$ is nonadjacent to $b_1$, but possibly $r_1=a$.
Let $R_1$ be the subpath of $P_1$ between $r_1,a$, and let $S_1$ be the subpath between $s_1,b_1$. There is a pyramid with apex $q$
and constituent paths $q\dd s_1\dd S_1\dd b_1$, $q\dd r_1\dd R_1\dd a\dd P_3$, and $P_2'$. Some two of these paths have lengths summing
to at least $\ell_1+\ell_2$. But the first two sum to $\ell_1+\ell_3<\ell_1+\ell_2$; the first and third sum to
at most $(E(S_1)|+1)+(\ell_2 -\ell_3-1)<\ell_1+\ell_2$; and the second and third sum to at most
$$(1+|E(R_1)|+\ell_3) + (\ell_2 -\ell_3-1)<\ell_1+\ell_2,$$
a contradiction. This proves (2).

\bigskip
From (1) and (2), it follows that $q$ has exactly two neighbours in $V(P_1)$ and they are adjacent. Let them be $r_1,s_1$,
where $a,r_1,s_1,b_1$
are in order in $P_1$. Let $R_1$ be the subpath of $P_1$ between $r_1,a$, and let $S_1$ be the subpath between $s_1,b_1$.
Let $P_2''$ be an induced path between $a, q$ with interior in $V(R_2\cup C_2\cup Q)$; and hence with length at most
$$\lceil \ell_2/2 \rceil + |E(Q)|\le (\lceil \ell_2/2 \rceil )+(\lceil \ell_2/2\rceil -\ell_3 -1)\le \ell_2-\ell_3.$$
There is a pyramid with apex $a$ and constituent paths $R_1$, $a\dd P_3\dd b_3\dd b_1\dd S_1 \dd s_1$, and
$P_2''$; and some two of them have lengths that sum to at least $\ell_1+\ell_2$. But the first two lengths sum to
$\ell_1+\ell_3<\ell_1+\ell_2$; the first and third sum to at most
$$|E(R_1)|+ (\ell_2-\ell_3)<\ell_1+\ell_2;$$
and the second and third sum to at most 
$$(\ell_3+1+|E(S_1)|)+ (\ell_2-\ell_3)<\ell_1+\ell_2$$
since $|E(S_1)|\le \ell_1-2$ (because $q,a$ are nonadjacent). This is a contradiction, and so proves \ref{chooseC2D2}.~\bbox

\section{Locating a great pyramid}

If $X\subseteq V(G)$, we define $N[X]$ to be the set of all vertices that either belong to $X$ or have a neighbour in $X$, and 
$N(X)=N[X]\setminus X$. If $X=\{v\}$ we write $N[X]$ for $N[\{x\}]$ and so on.
Now we are ready to prove \ref{findgreatpyr}, which we restate:
\begin{thm}\label{findgreatpyragain}
There is an algorithm with the following specifications:
\begin{description}
\item [Input:] A graph $G$.
\item [Output:] Outputs either an odd hole of $G$, or a statement of failure. If $G$ contains no 5-hole, and no jewelled shortest odd
hole, and $G$ contains a great pyramid,
the output will
be a shortest odd hole of $G$.
\item [Running time:] $O(|G|^{14})$.
\end{description}
\end{thm}
\Proof We enumerate all 12-tuples $(a,b_1,b_2,b_3, c_2,d_2,m_2,v,v_1,v_2,v_3,v_4)$ of vertices of $G$
such that $a,b_1,b_2,b_3$ are distinct and $b_1,b_2,b_3$ are pairwise adjacent.
For each one we carry out the following steps:
\begin{itemize}
\item Let $Y$ be $N[b_1]\cup (N[\{v,v_1,v_2\}]\setminus \{v_1,v_2,v_3,v_4\})$.
\item Let $X_1= Y\cup N[b_2]$. 
Choose a shortest path $Q_3$ between $a,b_3$ with interior in $V(G)\setminus X_1$. (If there is no such path, move on the next
$12$-tuple.)
\item Let $X_2=Y\cup N[V(Q_3)\setminus \{a\}]$. 
Choose a shortest path $R_3$ between $a,c_2$ with interior 
in $V(G)\setminus X_2$, and a shortest path $S_2$ between $b_2,d_2$ with interior 
in $V(G)\setminus X_2$. (If there are no such paths, move on.)
\item Let $X_3=Y\cup N[V(P_3)]$. Choose a 
a shortest path $R_2$ between $r_2,m_2$ with interior
in $V(G)\setminus X_4$, and a shortest path $S_2$ between $s_2,m_2$ with interior
in $V(G)\setminus X_4$ (and if there are no such paths, move on).
\item Let $X_4=N[V(R_2\cup S_2\cup C_2\cup D_2\cup P_3)\setminus \{a\}]$,
and choose a shortest path $Q_1$ between $a,b_1$ with interior
in $V(G)\setminus X_4$ (if there is no such path, move on).
\item Check whether $b_1b_2$ is an edge and the union of $P_1,R_2,C_2,D_2,S_2$ and the edge $b_1b_2$ is an odd hole, and if so, record it.
\end{itemize}
After examining all 12-tuples, if no hole is recorded we report failure, and otherwise output the shortest of the recorded holes.

To prove correctness, since the output is either failure or an odd hole, we only need check that when 
$G$ contains no 5-hole, and no jewelled shortest odd
hole, and $G$ contains a great pyramid, then 
the output will
be a shortest odd hole of $G$.
Thus, let $H$ be an optimal great pyramid, with apex $a$ and base $\{b_1,b_2,b_3\}$. Let its constituent paths be 
$P_1,P_2,P_3$ in the usual notation. Let $m_2\in V(P_2)$ such that the subpath of $P_2$ between 
$m_2,a$ has length $\lceil \ell_2/2\rceil$. If $\ell_2\ge 2\ell_3$, let $c_2\in V(P_2)$ such that the subpath of $P_2$ between
$a,c_2$ has length $\ell_3$, and define $d_2$ similarly; and otherwise let $c_2=d_2=m_2$.

We claim that there is a five-tuple $(v,v_1,v_2,v_3,v_4)$ of vertices such that, if we define $Y$ to be as in the first bullet above,
then every major vertex belongs to $Y$, and no vertex of $P_1^*\cup P_2^*$ belongs to $Y$. To see this, first suppose that every 
big $C$-major vertex belongs to $N[b_1]$;
then when we take $v=v_1=v_2=v_3=v_4=b_1$, the set $Y=N[b_1]$, and the claim holds. So we may assume that there
is a major vertex $v$ not in $N[b_1]$. Choose $v$ adjacent to $a$ if possible.
By \ref{pyrmajor}, and exchanging $P_1,P_2$ if necessary, we may assume that $v$ has type $(P_1,P_2)$
or $(P_1,P_3)$. In either case $v$ has exactly two neighbours in $V(P_1\cup P_2)$, say $p,q$. Also, by \ref{3.4mod}, there is
an edge $v_1v_2$ of $C$ such that $v$ is adjacent to one of $v_1,v_2$, and every other big $C$-major vertex not in $N[b_1]$ 
is adjacent to 
one of $v,v_1,v_2$; and therefore we may choose $v_1v_2$ with $v_1,v_2\ne b_1$. 
Choose $v_1v_2$ with $v_1,v_2\ne a$ if possible.

We claim that if one of $v_1,v_2=a$, then $v$ is adjacent to $a$. Suppose that $v_1=a$ say. If no big $C$-major vertex
is adjacent to $a$, then we can replace $v_1v_2$ by the other edge of $C$ that contains $v_2$, a contradiction. So some big $C$-major vertex is adjacent to $a$, and hence so is $v$, from the choice of $v$.
\\
\\
(1) {\em At most two vertices
of $(P_2^*\cup P_3^*)\setminus \{v_1,v_2\}$ are equal or adjacent to a member of $\{v,v_1,v_2\}$.}
\\
\\
To see this, there are several cases, depending on the position in $C$ of the edge $v_1v_2$. Let $Z$ be
the set of vertices
of $(P_2^*\cup P_3^*)\setminus \{v_1,v_2\}$ are equal or adjacent to a member of $\{v,v_1,v_2\}$. If $v_1,v_2\in V(P_1)\setminus \{a\}$,
then $Z=\{p,q\}$. If $v_1=a$ and $v_2\in V(P_1)$ or vice versa, then $v$ is adjacent to $a$, as we saw earlier, and $Z$ is the
set of the two neighbours of $a$ in $P_1\cup P_2$. If $v_1,v_2\in V(P_2)\setminus \{a\}$, then one of $p,q$ equals one of $v_1,v_2$
(since $v$ is adjacent to one of $v_1,v_2$) and the other of $p,q$ is adjacent in $P_2$ to one of $v_1,v_2$; so $Z$ is the set of 
the at most two vertices in $P_2$ that are adjacent to one of $v_1,v_2$ and different from both $v_1,v_2$. If $v_1=a$ and $v_2\in V(P_2)$,
then as before $v$ is adjacent to $a$ and therefore $\{p,q\}=\{v_1,v_2\}$, and $Z$ consists of the (at most two) vertices 
of $P_1\cup P_2)$ that are adjacent to one of $v_1,v_2$ and different from them both. The last case, when
$v_1=b_1$ and $v_2=b_2$, does not occur, because we chose $v_1,v_2\ne b_1$. This proves (1).

\bigskip

From (1), this proves that there is a five-tuple $(v,v_1,v_2,v_3,v_4)$ of vertices such that, 
when we define $Y$ to be as in the first bullet above,
every major vertex belongs to $Y$, and no vertex of $P_1^*\cup P_2^*$ belongs to $Y$.

We claim that when the algorithm examines this 12-tuple $(a,b_1,b_2,b_3, c_2,d_2,m_2,v,v_1,v_2,v_3,v_4)$, 
it will record a shortest odd hole.
To see this, let $Q_3$ be the path chosen in the first bullet above (it exists, since $P_3$ exists); then 
by \ref{chooseP3} we can replace $P_3$ by $Q_3$ and obtain another optimal great pyramid; that is, we can choose $H$ such that
$P_3=Q_3$.
By \ref{chooseC2D2}, the union of the four paths $R_2,C_2,D_2,S_2$
chosen in the second and third bullets above
is a path $Q_2$ with the same length as $P_2$, and we can choose $H$ such that $P_2=Q_2$ (while maintaining that $P_3=Q_3$). 
Now let $Q_1$
be as chosen in the fourth bullet (it exists, since $P_1$ exists); then it has length at most that of $P_1$, and forms a pyramid 
$H''$ with $P_2,P_3$; and $H''$ is a great pyramid, since $Q_1$ has length at most that of $P_1$ and $H$ is a great pyramid.
In particular $Q_1$ has the same length as $P_1$, and $G[V(Q_1\cup Q_2)]$ is a shortest odd hole, that the algorithm will record.
This proves correctness. For each $12$-tuple, the running time is $O(|G|^2)$, so the total running time is as claimed. This proves
\ref{findgreatpyr}.~\bbox

\section{Acknowledgement}
This work was done during the Structural Graph Theory Downunder Matrix Program, that was held
in 2019 at the Creswick Campus of the University of Melbourne. The authors express their gratitude
to the Matrix Center for the funding it provided and the use of its facilities.

\end{document}